\documentclass{article}
\usepackage{amsmath,amsfonts,amssymb,amsthm}
\usepackage{color}

\usepackage{hyperref}

\providecommand{\R}{\mathbb{R}}

\renewcommand{\leq}{\leqslant}
\renewcommand{\geq}{\geqslant}

\renewcommand{\div}{\operatorname{div}}
\newcommand{\curl}{\operatorname{curl}}

\newtheorem{Theorem}{Theorem}
\newtheorem{Definition}{Definition}

\newtheorem{Proposition}{Proposition}

\newtheorem{Remark}{Remark}

\begin{document}

\date{\today}
\title{A Kato type Theorem for the inviscid limit of the Navier-Stokes equations with a moving rigid body.  }

\author{Franck Sueur\footnote{Laboratoire Jacques-Louis Lions,
Universit\'e Pierre et Marie Curie - Paris 6,
4 place Jussieu,
75005 Paris, 
France}
}

\maketitle

\begin{abstract}

The issue of the inviscid limit for the incompressible Navier-Stokes equations  when a no-slip condition is prescribed on the boundary
is a famous open problem.
A   result  by Kato \cite{Tosio}  says that  convergence to the  Euler equations
holds  true in the energy space if and only if the energy dissipation
rate of the viscous flow  in a  boundary layer of width proportional to the viscosity vanishes. 
Of course, if one considers the motion of a solid body in an incompressible fluid, with a no-slip condition at the interface, the issue of the inviscid limit is as least as difficult.
However  it is not clear if the additional difficulties linked to the body's dynamic make this issue more difficult or not.
In this paper we consider the motion of a rigid body in an incompressible fluid occupying the complementary set in the space and we prove that a Kato type condition implies the convergence of the fluid velocity and of the body velocity as well, what seems to indicate that an answer in the case of a fixed boundary could also bring an answer to the case where there is a moving body in the fluid.

\end{abstract}

\section{Introduction}

In this paper we investigate the issue of the inviscid limit for a incompressible fluid, driven by 
the Navier-Stokes equations, in the case where there is a moving body in the fluid.
When a no-slip condition is prescribed on a solid boundary this issue is still widely open, even if this boundary does not move (see for instance \cite{constantin,BardosTiti,E,grenier}). 
However in this case a result by Kato \cite{Tosio} says that, in the inviscid limit,  the convergence to the  Euler equations
holds  true in the energy space if and only if the energy dissipation rate of the viscous flows  in a  boundary layer of width proportional to the viscosity vanishes. 
The main result in this paper is an extension of Kato's result in the case where there is a moving body in the fluid.
In order to clarify the presentation of our result we first recall Kato's result in its original setting: the case of a fluid contained in a fixed bounded domain, 
along with a slight reformulation which will be natural in the case with a moving body.

\subsection{A short review of Kato's result.}
Let us first consider the case of a fluid alone, contained in a bounded domain $\Omega \subset \R^d$, with $d=2$ or $3$. 
We therefore consider the incompressible  Navier-Stokes equations:
\begin{gather}
\label{NS1alone}
\frac{\partial U}{\partial t}+(U\cdot\nabla)U + \nabla P = \nu \Delta U \ \text{ for } \ x \in  \Omega, \\
\label{NS2alone}
\div U = 0 \ \text{ for } \  x \in \Omega  ,  \\
\label{NS3alone}
U = 0  \ \text{ for } \ x\in \partial  \Omega  ,  \\
\label{NS4alone}
U |_{t= 0} = U_0 .
\end{gather}
Here  $U$ and  $P$  denote respectively the velocity and pressure fields. The positive constant $\nu$ is the viscosity of the fluid.
The condition \eqref{NS3alone} is the so-called no-slip condition.

We are going to deal with weak solutions of \eqref{NS1alone}-\eqref{NS4alone}. Let us recall the following result by Leray (cf. for instance \cite{lions}), where we denote
\begin{eqnarray*}
 \mathcal{H}_\Omega &:=& \{ V \in L^{2} (\Omega) / \ \div V = 0    \text{ in } \Omega \text{ and } V \cdot n = 0  \text{ on }   \partial  \Omega\} , 
 \\  \mathcal{V}_\Omega &:=&   \{  V \in H^{1}_0 (\Omega) / \ \div V = 0  \text{ in }   \Omega  \} . 
\end{eqnarray*}

Let us warn here the reader that we  use the following slight abuse of notations: 
if $V$ denotes any scalar-valued function space and $U$ is a function with its values in $\R^d$, we will say that $U\in V$ if its components are  in $V$.

\begin{Theorem}
\label{leray2d}
Let $U_0 \in \mathcal{H}_\Omega$ and $T >0$. Then there exists a solution $U \in C_w (  [0,T ]  ; \mathcal{H}_\Omega ) \cap L^2 ( [0,T ] ;  \mathcal{V}_\Omega ) $ of the equations \eqref{NS1alone}-\eqref{NS4alone} in the sense that for all $V \in   H^1 ( [0,T ]  ; \mathcal{H}_\Omega ) \cap L^2 ( [0,T ] ;  \mathcal{V}_\Omega )$, for all $t\in [0,T ] $, 
\begin{eqnarray}
\label{WeakNSalone}
 \int_{\Omega } \Big( U (t,\cdot )  \cdot V  (t,\cdot ) -  U_{0} \cdot V |_{t=0} \Big) dx    = 
 \int_{0}^{t} \int_{\Omega } \Big[ U  \cdot  (\partial_{t} +U   \cdot\nabla) V      - 2 \nu  \nabla U :  \nabla V  \Big]  dx  ds .
\end{eqnarray}
Moreover this solution satisfies the following energy inequality: for any $t \in [0,T ]$, 
\begin{eqnarray}
\label{NSBodyWeakEnergyAlone}
\frac{1}{2} \| U (t, \cdot)  \|^{2}_{L^{2} (\Omega) } + \nu \int_{(0,t) \times \Omega } |   \nabla U |^{2} dx ds \leq \frac{1}{2}  \| U_{0}  \|^{2}_{L^{2} (\Omega)} .
\end{eqnarray}
Moreover when $d=2$ this solution is unique, $U \in C (  [0,T ]  ; \mathcal{H}_\Omega )$ and there is equality in \eqref{NSBodyWeakEnergyAlone}.
\end{Theorem}
When the viscosity coefficient $\nu$ is set equal to $0$ in the previous equations, it is expected that the system \eqref{NS1alone}-\eqref{NS4alone} degenerates into the following incompressible  Euler equations:
\begin{gather}
\label{Euler1alone}
\frac{\partial U^E}{\partial t} + (U^E \cdot\nabla) U^E  + \nabla P^E  =  0  \ \text{ for } \ x \in \Omega , \\
\label{Euler2alone}
\div U^E  = 0 \ \text{ for } \  x \in  \Omega,  \\
\label{Euler3alone}
U^E  \cdot n =  0  \ \text{ for } \ x\in \partial  \Omega ,  \\
\label{Euler4alone}
U^E  |_{t= 0} = U^E_0 .
\end{gather}
Kato's result deals with classical solutions of the Euler equations  \eqref{Euler1alone}-\eqref{Euler4alone}, whose (local in time)  existence and uniqueness are classical since the works of  Lichtenstein, G\"{u}nter and Wolibner.
Let us also recall that in two dimensions they are global in time, cf.   \cite{wolibner} in the case of a simply connected domain and  \cite{kato} for  multiply connected domains.
More precisely we have the following result, where we make use of the notation  $C^{1,\lambda}(\Omega)$ for the H\"older space, endowed with the norm:
\begin{align*}
\| V  \|_{  C^{1,\lambda} ( \Omega  ) } := 
   \|   V  \|_{L^\infty (  \Omega ) }  
+ \sup_{ x \neq y \in   \Omega  } \frac{ |\nabla  V (x) -  \nabla V (y)| }{  |x - y|^\lambda }   .
\end{align*}
Here $\lambda \in (0,1)$.
\begin{Theorem}
\label{Strongalone}
Let be given $U^E_0 \in \mathcal{H}_\Omega \cap C^{1,\lambda}(\Omega)$.
Then there exists $T >0$ and a unique solution $U^E $ of \eqref{Euler1alone}-\eqref{Euler4alone} in $  C ( [0,T ]  ;  \mathcal{H}_\Omega) \cap  C_{w*} ( [0,T ]  ; C^{ 1,\lambda}(\Omega ))$.
Moreover this solution satisfies the following energy equality: for any $t \in [0,T ]$, 
\begin{eqnarray}
\label{EulerBodyStrongEnergyAlone}
\| U^E  (t, \cdot)  \|_{ L^2 (  \Omega )} = \|  U^E _{0}  \|_{L^2 (  \Omega )} .
\end{eqnarray}
Moreover in two dimensions, $T $ can be chosen arbitrarily.
\end{Theorem}

We are now in position to recall Kato's result.
\begin{Theorem}
\label{KatoAlone}
Let be given $c>0$ and $T >0$. 
Assume that 
$U^E_{0} \in  \mathcal{H}_\Omega \cap  C^{ 1,\lambda}(\Omega)$ and that 
$U_{0} \rightarrow  U^E_{0} \text{ in }  \mathcal{H}_\Omega$ when  $\nu \rightarrow  0$.  
Let us denote by $U$ a solution  of \eqref{NS1alone}-\eqref{NS4alone} given by Theorem \ref{leray2d} and by  $ U^E $ the solution of \eqref{Euler1alone}-\eqref{Euler4alone} given by Theorem  \ref{Strongalone}.
Let us denote
\begin{eqnarray*}
\Gamma^{\Omega}_{c\nu} := \{  x \in  \Omega / \  dist (x,  \partial  \Omega  )  < c \nu \} ,
\end{eqnarray*}
which is well defined for $\nu >0$ small enough.

Then the following conditions are equivalent, when  $\nu \rightarrow  0$.
\begin{enumerate}
\item $\nu \int_{(0,T) \times \Gamma^{\Omega}_{c\nu} } | \nabla U |^{2} dx dt \rightarrow  0 $,
\item  $U \rightarrow U^E $   in  $  C ( [0,T ]  ;   \mathcal{H}_\Omega )$. 
\end{enumerate}
\end{Theorem}
Comparing \eqref{NSBodyWeakEnergyAlone} and \eqref{EulerBodyStrongEnergyAlone} we see that the quantity in the first condition in Theorem \ref{KatoAlone} can be interpreted as the energy dissipation rate of the viscous flows in a boundary layer of width proportional to the viscosity.
This width is much smaller than the one given by Prandtl's theory, what seems to indicate that one has to go beyond Prandtl's description to understand 
 the inviscid limit. Moreover some recent results 
 \cite{DGVD,guo} show that Prandtl's equation is in general ill-posed.

Kato's result in \cite{Tosio} contains some extra considerations about source terms and weak convergence, but we will skip these considerations here for sake of simplicity.
Furthermore there exists many variants of Kato's argument: see for instance \cite{Wang,TemamWang,Kelliher,Masmoudi,ILL}. 
In particular it is shown in  \cite{Kelliher} that another equivalent condition is\footnote{Here we use the following notations: when $A$ and $B$ are two $d \times d$  matrices,  we denote $A:B = \sum_{1 \leqslant i,j \leqslant  d} A_{ij}B_{ij}$ and $| A |^2:= A:A $.}
$$\nu \int_{(0,T) \times \Gamma^{\Omega}_{c\nu} } | \curl U |^{2} dx dt \rightarrow  0 ,$$
where $ \curl U$ is the $d \times d$ skew symmetric matrix given by
$$  \curl U := ( \frac{1}{2} ( \partial_{j} U_{i} -  \partial_{i} U_{j} ) )_{1 \leqslant i,j \leqslant  d} ,$$
and a slight modification of the proof  in  \cite{Kelliher}  also yields that another  equivalent condition is 
\begin{eqnarray}
\label{defo2}
\nu \int_{(0,T) \times \Gamma^{\Omega}_{c\nu} } | D(U) |^{2} dx dt \rightarrow  0 ,
\end{eqnarray}
where $D(U) $ is the deformation tensor
\begin{eqnarray}
\label{defo}
 D(U) := ( \frac{1}{2} ( \partial_{j} U_{i} +  \partial_{i} U_{j} ) )_{1 \leqslant i,j \leqslant  d} .
\end{eqnarray}
Actually, the proof in \cite{Kelliher} relies on the observations that 
\begin{enumerate}
\item for any $U,V$ in $H^1 (\Omega )$ such that $\div V = 0$  and such that $ ( U \cdot \nabla V) \cdot n = 0$ on $ \partial  \Omega  $,
\begin{eqnarray}
\label{P1curl}
\int_{\Omega } \nabla U : \nabla V = 2 \int_{\Omega }  \curl U  :  \curl V ,
\end{eqnarray}
\item  for any $U$ in $H^1 (\Omega )$ and $V$ in $C^1  (\Omega )$ such that $\div V = 0$ and such that $(n \cdot V) U = 0$ on $ \partial  \Omega  $,
\begin{eqnarray}
\label{P1grad}
\int_{\Omega } V \cdot ( U \cdot \nabla U )   = 2 \int_{\Omega } V \cdot \Big( (\curl U ) U \Big) .
\end{eqnarray}
\end{enumerate}
 These properties also hold true when we substitute $ D(U) $ to $ \curl U$, that is 
\begin{enumerate}
\item for any $U,V$ in $H^1 (\Omega )$ such that $\div V = 0$  and such that $ ( U \cdot \nabla V) \cdot n = 0$ on $ \partial  \Omega  $,
\begin{eqnarray}
\label{P1}
\int_{\Omega } \nabla U : \nabla V = 2 \int_{\Omega }  D(U)  :  D(V) ,
\end{eqnarray}
\item  for any $U$ in $H^1 (\Omega )$ and $V$ in $C^1  (\Omega )$ such that $\div V = 0$ and such that $(n \cdot V) U = 0$ on $ \partial  \Omega  $,
\begin{eqnarray}
\label{P2}
\int_{\Omega } V \cdot ( U \cdot \nabla U )   = 2 \int_{\Omega } V \cdot \Big( D(U) U \Big) .
\end{eqnarray}
\end{enumerate}

 It is therefore sufficient to follow the proof in \cite{Kelliher} with these subtitutions in order to add  \eqref{defo2} to the list of the equivalent conditions in Theorem \ref{KatoAlone}.
We are going to use a condition similar to \eqref{defo2} in the case of a moving rigid body.

\subsection{The case of a fluid with a moving rigid body.}

We now consider the case where there is a moving rigid body in a fluid. Let us focus here on the three dimensional case. 
We assume that the body initially occupies  a closed, bounded, connected and simply connected subset $\mathcal{S}_0 \subset \R^3$ with smooth boundary.
It rigidly moves so that at  time $t$  it occupies an isometric  domain denoted by $\mathcal{S}(t)$.
More precisely if we denote by $h (t)$ the position of the center of mass of the body at time $t$, then 
there exists a rotation matrix $Q (t) \in SO(3)$, 
such that  the position $\eta (t,x) \in \mathcal{S} (t)$  at 
the time $t$ of the point fixed to the body with an initial position $x$ is 
\begin{eqnarray} \label{FlotSolide}
\eta (t,x) := h (t) + Q (t)(x- h (0)) .
\end{eqnarray}
Of course this yields that $Q(0) =0$.
Since  $Q^{T} Q' (t) $ is skew symmetric there exists (only one) $r  (t)$ in $\R^3$ such that for any $x \in \R^3$, 
\begin{equation} \label{LoiDeQ}
Q^{T} Q' (t)  x = r(t) \wedge x .
\end{equation}
Accordingly, the solid velocity is given by
\begin{equation*} 
U_{{\mathcal S}}(t,x) := h'(t)  + R(t) \wedge (x-h(t)) \text{ with } R(t) := Q(t)r(t) .
\end{equation*}
Given a positive function $\rho_{S_{0}} $, say in $ L^{\infty}({\mathcal S}_{0};\R)$, describing the density in the solid, the solid mass $m>0$, the center of mass $h(t)$ and the  inertia matrix ${\mathcal J} (t)$ can be computed by it first moments. Let us recall that $ {\mathcal J} (t)$ is symmetric positive definite and that  $\mathcal{J}$ satisfies  Sylvester's law:
\begin{equation}
\label{Sylvester}
 \mathcal{J}(t) = Q(t)  \mathcal{J}_0 Q^{T} (t) ,
 \end{equation}
where $\mathcal{J}_0$ is the initial value of $\mathcal{J}$.

 In the rest of the plane, that is in the open set 
$\mathcal{F}(t) := \R^3 \setminus {\mathcal S} (t)$,
evolves a planar ideal fluid driven by the incompressible Navier-Stokes equations. 
We denote correspondingly ${\mathcal F}_{0}:=\R^{3} \setminus {\mathcal S}_{0}$ the initial fluid domain. \par
The complete system driving the dynamics reads
\begin{gather}
\label{NS1}
\frac{\partial U}{\partial t}+(U\cdot\nabla)U + \nabla P = \nu \Delta U + g \ \text{ for } \ x \in \mathcal{F}(t), \\
\label{NS2}
\div U = 0 \ \text{ for } \  x \in \mathcal{F}(t) ,  \\
\label{NS3}
U =  U_\mathcal{S}  \ \text{ for } \ x\in \partial \mathcal{S}  (t),  \\
\label{Solide1} 
m  h'' (t) = m g  - \int_{ \partial \mathcal{S} (t)} \Sigma n \, ds ,  \\
\label{Solide2} 
(\mathcal{J}  R )' (t) = -    \int_{ \partial   \mathcal{S} (t)}  (x-  h )  \wedge \Sigma  n \, ds , \\
\label{NSci2}
U |_{t= 0} = U_0 , \\
\label{Solideci}
h (0)=  0 , \ h' (0)= \ell_0, \ R  (0)=  r_0. 
\end{gather}
Here  $U$ and  $P$  denote the fluid velocity and pressure, which are   defined on ${\mathcal F}(t)$ for each $t$.
  The fluid is supposed to be  homogeneous of density $1$, to simplify the notations and without any loss of generality. 
The Cauchy stress tensor is defined by
\begin{eqnarray*}
 \Sigma := -P Id + 2 \nu D(U) ,
\end{eqnarray*}
where $D(U)$ is the deformation tensor defined in \eqref{defo}.
\par
Above $n$ denotes the unit outward normal on the boundary of the fluid domain, 
$ds$ denotes the integration element on this boundary and $g$ is the gravity force which is assumed to be a constant vector, we actually include it in our study as a physical  example of source term.

Let us observe that the choice $h (0)=0$ avoids to write an extra moment, the one due to the gravity force, in \eqref{Solide2}.
Still this choice is only a matter of convention and does not decrease the generality.

The existence of a weak solution to the system \eqref{NS1}-\eqref{Solideci} was given in \cite{Serre}.
Let us also refer here to the following subsequent works  \cite{DE1,DE2,Conca1,Conca2,TakTuc} and the references therein. \par
\par

When the viscosity coefficient $\nu$ is set equal to $0$ in the previous equations, formally, the system \eqref{NS1}-\eqref{Solideci} degenerates into the following equations:
\begin{gather}
\label{Euler1}
\frac{\partial U^E}{\partial t} + (U^E \cdot\nabla) U^E  + \nabla P^E  =   g   \ \text{ for } \ x \in \mathcal{F}^E  (t), \\
\label{Euler2}
\div U^E  = 0 \ \text{ for } \  x \in \mathcal{F}^E  (t) ,  \\
\label{Euler3}
U^E  \cdot n =  U^E_{{\mathcal S}}   \cdot n  \ \text{ for } \ x\in \partial \mathcal{S}^E   (t),  \\
\label{EulerSolide1} 
m   (h^E) ''  =  m g + \int_{ \partial \mathcal{S}^E  (t)}  P^E  n \, ds ,  \\
\label{EulerSolide2} 
(\mathcal{J}^E R^E   )' =    \int_{ \partial   \mathcal{S}^E (t)}  P^E  (x-  h^E ) \wedge n \, ds    , \\
\label{Eulerci2}
U^E  |_{t= 0} = U^E_0 , \\
\label{EulerSolideci}
h^E  (0)=  0 , \ (h^E )' (0)= \ell^{E}_0 ,\ R^{E}  (0)=  r^{E}_0 ,
\end{gather}
where
\begin{equation*} 
U^E_{{\mathcal S}}(t,x) := (h^E) '(t)  +  R^E (t)  \wedge  (x-h^E(t)) ,
\end{equation*}
and
\begin{equation*} 
\mathcal{S}^E   (t) := \eta^E (t,\cdot) (\mathcal{S}_0 ) ,\text{ with } \eta^E (t,x) := h^E (t) + Q^E (t) x ,
\end{equation*}
where the matrix $ Q^E$ solves the  differential equation 
$ (Q^E)'    = R^E  \wedge Q^E $ with $Q^E (0) = 0$. Finally $ \mathcal{J}^E$ is given by
$ \mathcal{J}^E = Q^E \mathcal{J}_0 (Q^E)^{T} $.

Observe that we prescribe $h^E  (0)=  0$ so that the initial position $ \mathcal{S}^E (0)$ occupied by the solid also starts from $ \mathcal{S}_0$ at $t=0$.
The mass $m$ and the initial inertia matrix $ \mathcal{J}_0 $ are also the same than in the previous case of the Navier-Stokes equations.

The existence and uniqueness of classical solutions to the equations \eqref{Euler1}-\eqref{EulerSolideci} is now well understood thanks to the recent works \cite{Ortega,Ortega2,rosier,ht,ogfstt,GS}.\\
\\
The aim of this paper is to show the following conditional result about the inviscid limit:  if 
\begin{eqnarray}
\label{dissiP}
 \nu \int_{(0,T)}  \int_{  \Gamma_{c\nu} (t) } | D(U) |^{2} dx dt  \rightarrow  0 ,
\end{eqnarray}
when  $\nu \rightarrow  0$, where, for some $c>0$, 
\begin{eqnarray*}
\Gamma_{c\nu} (t):= \{  x \in  \mathcal{F} (t) / \  dist (x,  \mathcal{S} (t) )  < c \nu \}  ,
\end{eqnarray*}
then the solution  of \eqref{NS1}-\eqref{Solideci} converges  to the solution  of \eqref{Euler1}-\eqref{EulerSolideci}.

A precise statement is given below. In particular we will see that the condition  \eqref{dissiP} is also necessary.

\section{Change of variables}
In order to write the equations of the fluid in a fixed domain, we are going to use some changes of variables.
\subsection{Case of the Navier-Stokes equations}
In the case of the Navier-Stokes equations we use the following change of variables:
\begin{eqnarray*}
\ell(t) :=   Q(t)^T \,   h' (t)  ,  \, u(t,x) := Q(t)^T \,  U(t, Q(t) x+h(t)) ,
\\ p(t,x) := P(t, Q(t) x+ h(t)) 
 \text{ and } \sigma (t,x) := \Sigma (t, Q(t)  x+ h(t)) ,
\end{eqnarray*}
so that
\begin{eqnarray*}
 \sigma := -p Id + 2 \nu D(u) , \text{ where  }  D(u) := ( \frac{1}{2} ( \partial_{j} u_{i} +  \partial_{i} u_{j} ) )_{i,j} .
\end{eqnarray*}
Therefore the system \eqref{NS1}-\eqref{Solideci}  now reads 
\begin{gather}
\label{chNS1}
\frac{\partial u}{\partial t}+ (  u - u_\mathcal{S} ) \cdot \nabla u +  r \wedge u  + \nabla p = Q(t)^T \,  g  + \nu \Delta u \ \text{ for } \ x \in \mathcal{F}_0 , \\
\label{chNS2}
\div u = 0 \ \text{ for } \  x \in \mathcal{F}_0  ,  \\
\label{chNS3}
u  = u_\mathcal{S} \ \text{ for } \ x\in \partial \mathcal{S}_0,  \\
\label{chSolide1} 
m  \ell'  = m  Q^T \,  g  - \int_{ \partial \mathcal{S}_0 } \sigma n \, ds  + m \ell \wedge  r,  \\
\label{chSolide2} 
\mathcal{J}_0 r'  =   - \int_{ \partial   \mathcal{S}_0}   x \wedge \sigma n \, ds  + ( \mathcal{J}_0  r) \wedge r , \\
\label{chNSci2}
u |_{t= 0} = u_0 , \\
\label{chSolideci}
h (0)= 0  , \ h' (0)= \ell_0 ,\ r  (0)=  r_0. 
\end{gather}
with
\begin{eqnarray}
\label{LastReading}
  u_{\mathcal{S}} (t,x) := \ell  (t)  + r(t) \wedge x .
   \end{eqnarray}
In order to write the weak formulation of the system \eqref{chNS1}-\eqref{chSolideci} we introduce
\begin{eqnarray*}
 \mathcal{H} := \{ \phi \in L^{2} (\R^{3}) / \ \div \phi = 0  \text{ in }  \R^{3} \text{ and } D(\phi) = 0  \text{ in }  \mathcal{S}_0 \} .
\end{eqnarray*}
According to Lemma $1.1$ in \cite{Temam}, p18,  for all  $\phi \in  \mathcal{H}$, there exists $\ell_{\phi} \in \R^{3}$ and $r_{\phi} \in \R^{3} $ such that for any $x \in  \mathcal{S}_0$,
$\phi (x) = \ell_{\phi} + r_{\phi} \wedge  x $.
Therefore we extend the initial data $u_{0}$ (respectively $u^E_{0}$)  by setting 
$u_{0} := \ell_{0} + r_{0} \wedge x $ (resp.  $u^E_{0} := \ell^E_{0} + r^E_{0} \wedge x $) for $x \in \mathcal{S}_0 $.

We endow the space $L^{2} (\R^{3})$ with the following inner product:
\begin{eqnarray*}
(\phi  , \psi  )_{\mathcal{H}} :=    \int_{\mathcal{F}_0 }  \phi   \cdot   \psi  dx   +  \int_{\mathcal{S}_0 }  \rho_{S_{0}}  \phi   \cdot   \psi  dx .
\end{eqnarray*}
When $\phi  , \psi $  are in  $\mathcal{H}$ then,
\begin{eqnarray*}
(\phi  , \psi  )_{\mathcal{H}} =    \int_{\mathcal{F}_0 }  v_{\phi} \cdot  v_{\psi} dx  +  m  l_{\phi} \cdot  l_{\psi} + \mathcal{J}_0  r_{\phi}   \cdot   r_{\psi}  ,
\end{eqnarray*}
by definition of $m$ and $\mathcal{J}_0$.

\begin{Proposition}
\label{WeakStrong}
A smooth solution of \eqref{chNS1}-\eqref{chSolideci} satisfies the following: for any $v \in C^{\infty } ([0,T ] ;  \mathcal{H} \cap C^{\infty }_{c} (\R^{3} ))$, 
for all $t\in [0,T ] $, 
\begin{eqnarray}
\label{WeakNS}
(u,v)_{\mathcal{H}}  (t) -  (u_{0},v |_{t=0})_{\mathcal{H}}   = 
 \int_{0}^{t} \Big[ (u, \partial_{t} v)_{\mathcal{H}} 
  +  b(u,u,v)    - 2 \nu   \int_{\mathcal{F}_0 }     D(u) :  D(v)  dx + f_s [u,v]  \Big] ds ,
\end{eqnarray}
with
\begin{eqnarray*}
 f_t [u,v] :=  m_a   Q(t)^T \,  g  \cdot  \ell_{v} - Vol (\mathcal{S}_0 )  Q(t)^T \,  g  \cdot  (r_{v} \wedge x_0 ) ,
\end{eqnarray*}
where
\begin{eqnarray*}
 m_a := m - Vol (\mathcal{S}_0 ) \text{ and } x_0 := (Vol(\mathcal{S}_0))^{-1} \,   \int_{ \mathcal{F}_0} x dx
\end{eqnarray*}
are respectively the apparent mass and the centro\"id  of the solid,
and
\begin{eqnarray*}
b(u,v,w) := m  \det (r_u ,\ell_{v} ,  \ell_w ) + \det (\mathcal{J}_0 r_u , r_{v} , r_w ) +
\int_{\mathcal{F}_0 } \Big(     [   (u-  u_\mathcal{S} )  \cdot\nabla w ] \cdot v   -  \det (r_u , v ,  w ) \Big) dx
\end{eqnarray*}
\end{Proposition}
Let us stress that $ f_t [u,v] $ depends on $u$ via the rotation matrix $Q(t)$ which is obtained by solving the matrix differential equation
\begin{eqnarray}
\label{exp1}
Q' = Q (r_u \wedge \cdot)  \text{ with } Q(0)=Id .
\end{eqnarray}

We postpone the proof of Proposition \ref{WeakStrong} to the Appendix.
For the sequel we will need to enlarge the space of the test functions. Therefore we introduce the space
\begin{eqnarray*}
 \mathcal{V } :=   \{  \phi \in   \mathcal{H} / \ \int_{ \R^{3} } | \nabla \phi  (y) |^2 (1 + | y |^2 ) dy < + \infty    \} . 
\end{eqnarray*}
It is worth to notice from now on that $b$ is well-defined and trilinear on $\mathcal{H} \times \mathcal{H} \times \mathcal{V}$ (the weight above allowing to handle the rotation part of $u_\mathcal{S}$). Moreover it satisfies the following crucial property 
\begin{eqnarray}
\label{Tonga}
(u,v) \in \mathcal{H} \times \mathcal{V} \text{ implies } b(u,v,v) = 0 .
\end{eqnarray}
\begin{Definition}[] 
\label{Weak}
We say that 
\begin{eqnarray*}
u \in C_w ( [0,T ]  ; \mathcal{H} ) \cap L^2 ( [0,T ] ;  H^1 (\R^2 ) ) 
\end{eqnarray*}
is a weak solution of the system \eqref{chNS1}-\eqref{chSolideci} 
if for all $v \in   H^1 ( [0,T ]  ; \mathcal{H} ) \cap L^2 ( [0,T ] ;  \mathcal{V} )$,  
and for all $t\in [0,T ] $, \eqref{WeakNS} holds true. 
\end{Definition}
As already said above the existence of weak solutions "\`a la Leray" for the system \eqref{chNS1}-\eqref{chSolideci} is now well understood. 
Let us for instance refer to \cite{Serre}, Theorem 4.5.
\begin{Theorem}
\label{NSBodyWeak}
Let be given ${u}_{0} \in  \mathcal{H}$  and $T >0$. 
Then there exists a weak solution ${u}$ of \eqref{chNS1}-\eqref{chSolideci} in $  C_w ( [0,T ]  ;  \mathcal{H}) \cap  L^2 ( [0,T ] ; H^1 (\R^2 ) )$.
Moreover this solution satisfies the following energy inequality: for any $t \in [0,T ]$, 
\begin{eqnarray}
\label{NSBodyWeakEnergy}
\frac{1}{2} \| u (t, \cdot)  \|^{2}_{\mathcal{H}} + 2 \nu \int_{(0,t) \times \R^{3} } | D(u) |^{2} dx dt 
\leq \frac{1}{2}   \| u_{0}  \|^{2}_{\mathcal{H}} +  \int_{0}^{t}   f_s [u,u]  ds    .
\end{eqnarray}
\end{Theorem}
Let us stress that the integral above could innocuously be taken over $(0,t) \times \mathcal{F}_0$ since the deformation tensor $D(u) $ vanishes in the solid.

\begin{Remark}
\label{NSBodyWeakRemarkVort}
In the previous statement, it is possible to replace the weak formulation  \eqref{WeakNS} by the following one, based on the vorticity:
 for any $v \in C^{\infty } ([0,T ] ;  \mathcal{H} \cap C^{\infty }_{c} (\R^{3} ))$, 
for all $t\in [0,T ] $, 
\begin{eqnarray}
\label{WeakNSvort}
(u,v)_{\mathcal{H}}  (t) -  (u_{0},v |_{t=0})_{\mathcal{H}}   = 
 \int_{0}^{t} \Big[ (u, \partial_{t} v)_{\mathcal{H}} 
  +  b(u,u,v)    -  2 \nu \int_{\mathcal{F}_0 }     \curl u :  \curl v  dx + f_s [u,v]  \Big] ds ,
\end{eqnarray}
and the energy inequality \eqref{NSBodyWeakEnergy} by
\begin{eqnarray}
\label{NSBodyWeakEnergyvort}
\frac{1}{2} \| u (t, \cdot)  \|^{2}_{\mathcal{H}} + 2 \nu  \int_{(0,t) \times  \mathcal{F}_0} |  \curl u  |^{2} dx dt \leq \frac{1}{2}   \| u_{0}  \|^{2}_{\mathcal{H}} +  \int_{0}^{t}   f_s [u,u]  ds  .
\end{eqnarray}
\end{Remark}
\begin{Remark}
\label{NSBodyWeakRemarkGrad}
In Theorem \eqref{NSBodyWeak}, it is also possible to replace  \eqref{WeakNS} by:
 for any $v \in C^{\infty } ([0,T ] ;  \mathcal{H} \cap C^{\infty }_{c} (\R^{3} ))$, 
for all $t\in [0,T ] $, 
\begin{eqnarray}
\label{WeakNSgrad}
(u,v)_{\mathcal{H}}  (t) -  (u_{0},v |_{t=0})_{\mathcal{H}}   = 
 \int_{0}^{t} \Big[ (u, \partial_{t} v)_{\mathcal{H}} 
  +  b(u,u,v)    -  \nu \int_{\mathcal{F}_0 }     \nabla u :  \nabla v  dx + f_s [u,v]  \Big] ds ,
\end{eqnarray}
and  \eqref{NSBodyWeakEnergy} by
\begin{eqnarray}
\label{NSBodyWeakEnergygrad}
\frac{1}{2} \| u (t, \cdot)  \|^{2}_{\mathcal{H}} +  \nu \int_{(0,t) \times \mathcal{F}_0 } |  \nabla u  |^{2} dx dt \leq \frac{1}{2}   \| u_{0}  \|^{2}_{\mathcal{H}} +  \int_{0}^{t}   f_s [u,u]  ds .
\end{eqnarray}
\end{Remark}

\subsection{Case of the Euler equations}

Let us now see the case of the Euler equations. First
performing the following change of variables:
\begin{eqnarray*}
\ell^E (t) :=   Q^E(t))^T \,   (h^E )' (t)  ,  \, R^E(t) := Q^E(t)r^E(t) ,  \, 
\\ u^E (t,x) := Q^E (t)^T \,  U^E(t, Q^E(t) x+h^E(t)) ,
\text{ and }  p^E(t,x) := P^E(t, Q^E(t) x+ h^E(t)) ,
\end{eqnarray*}
where $Q^E(t)$ is 
 the rotation matrix associated to the motion of ${\mathcal{S}}^E (t)$,
the system  \eqref{Euler1}-\eqref{EulerSolideci}   now reads  
\begin{gather}
\label{chEuler1}
\frac{\partial u^E }{\partial t} +  ( u^E  - u^E_\mathcal{S} ) \cdot\nabla  u^E  +  r^E  \wedge u^E  + \nabla p^E  = Q^E g \ \text{ for } \ x \in \mathcal{F}_0  , \\
\label{chEuler2}
\div u^E  = 0 \ \text{ for } \  x \in \mathcal{F}_0  ,  \\
\label{chEuler3}\
u^E (t,x)  \cdot n = u^E_\mathcal{S}   \cdot n  \ \text{ for } \ x\in \partial \mathcal{S}_0  ,  \\
\label{chEulerSolide1} 
m   (\ell^E)' =  \int_{ \partial \mathcal{S}_0 } p^E  n \, ds + (m\ell^E ) \wedge  r^E ,  \\
\label{chEulerSolide2} 
\mathcal{J}_0   (r^E)'   =   \int_{ \partial   \mathcal{S}_0} p^E x \wedge  
n \, ds  +    ( \mathcal{J}_0  r^E) \wedge r^E  , \\
\label{chEulerci2}
u^E  |_{t= 0} = u^E _0 , \\
\label{chEulerSolideci}
h^E  (0)= 0 , \ (h^E) ' (0)= \ell^E_0 , \ r^E  (0)=  r^E_0 ,
\end{gather}
with
\begin{eqnarray}
\label{LastReadingE}
  u^E_{\mathcal{S}} (t,x) := \ell^E  (t)  + r^E (t) x^{\perp} .
   \end{eqnarray}

Here, in order to follow Kato's strategy we will need classical solutions.
The existence and uniqueness of classical solutions to the equations \eqref{Euler1}-\eqref{EulerSolideci}  with finite energy is given by the following result. 
\begin{Theorem}
\label{EulerBodyStrong}
Let be given  $\lambda \in (0,1)$ and $u^E_{0} \in  \mathcal{H}$ such that $ u^E_{0} |_{{\mathcal{F}_0}} \in  H^{1} \cap C^{ 1,\lambda} $ and $\curl u^E_{0}|_{{\mathcal{F}_0}} $ is compactly supported.
Let  $T >0$.
Then there exists a unique solution $u^E $ of \eqref{chEuler1}-\eqref{chEulerSolideci} in $  C^{1} ( [0,T ]  ;  \mathcal{H})$ such that 
$(\nabla u^E ) |_{ [0,T ]  \times {\mathcal{F}_0}} \in    C ( [0,T ]  ;  L^2 ({\mathcal{F}_0} , (1+ | x |^2 )^\frac12  dx )) \cap  C_{w*} ( [0,T ]  ; C^{0,\lambda} ({\mathcal{F}_0}))$.
Moreover for any $t \in [0,T ]$, 
\begin{eqnarray}
\label{EulerBodyStrongEnergy}
\frac12 \| u^E  (t, \cdot)  \|^2_{\mathcal{H}} = \frac12 \|  u^E _{0}  \|^2_{\mathcal{H}} +  \int_{0}^{t}   f_s [u^E , u^E]  ds  ,
\end{eqnarray}
where we denote, for $s  \in [0,T ]$ and $v \in \mathcal{H}$,
\begin{eqnarray*}
 f_s [u^E , v] :=  m_a   Q^E (s)^T \,  g  \cdot  \ell_{v} - Vol (\mathcal{S}_0 )  Q^E(s)^T \,  g  \cdot  (r_{v} \wedge x_0 ) .
\end{eqnarray*}
where the rotation matrix $Q^E (t)$ is obtained by solving the matrix differential equation
\begin{eqnarray}
\label{exp2}
 (Q^E) ' = Q^E (r^E \wedge \cdot) \text{ with } Q^E (0)=Id .
\end{eqnarray}
\end{Theorem}
Theorem \ref{EulerBodyStrong} can be proved in the same way than Th. $4$ in \cite{ogfstt}. The only difference is that Th. \ref{EulerBodyStrong} deals with the case where the fluid-rigid body system occupies the whole space whereas it was assumed to occupy a bounded domain in \cite{ogfstt}. Let us therefore only briefly discuss the decreasing at infinity of the fluid velocity in Theorem \ref{EulerBodyStrong}. Since the vorticity is transported (and stretched) by the flow and assumed to be compactly supported initially, it is compactly supported at any time. Then the fluid velocity $u$ can be recovered from the vorticity  by a Biot-Savart type operator, so that $u$ decreases as $x^{-2}$ at infinity and $\nabla_{x} u$ decreases as $x^{-3}$, uniformly in time. This entails the desired decreasing properties for $u$.

\section{Statement of the main result}
Let us now state the main result of this paper.
\begin{Theorem}[] 
\label{KatoBody}
Let be given $c>0$, $T >0$ and $u^E_{0} $ as in Theorem \ref{EulerBodyStrong}. 
Assume that 
\begin{eqnarray}
 \label{cvDonnees}
u_{0} \rightarrow  u^E_{0} \text{ in }  \mathcal{H}  \text{ when  }  \nu \rightarrow  0.
\end{eqnarray}
Let us denote ${u}$ a solution of \eqref{chNS1}-\eqref{chSolideci} given by Theorem \ref{NSBodyWeak} and by  $u^E$ the solution  of \eqref{chEuler1}-\eqref{chEulerSolideci} given by Theorem  \ref{EulerBodyStrong}.

Let us introduce the strips 
\begin{eqnarray*}
\Gamma_{c\nu} := \{  x \in  \mathcal{F}_0 / \  d(x)  < c \nu \} \text{ with }  d(x) := dist (x,  \partial   \mathcal{S}_0 ) ,
\end{eqnarray*}
which are well-defined for $\nu $ small enough.

Then the following conditions are equivalent, when  $\nu \rightarrow  0$:
\begin{eqnarray}
\label{ConV}
u \rightarrow u^E  \text{ in }   C ( [0,T ]  ;   \mathcal{H} ) ,
\\  \label{KatoCondition}
\nu \int_{(0,T) \times \Gamma_{c\nu} } | D(u) |^{2} dx dt \rightarrow  0 ,
\\ \label{KatoConditionVort}
\nu \int_{(0,T) \times \Gamma_{c\nu} } | \curl u |^{2} dx dt \rightarrow  0 ,
\\ \label{KatoConditionGrad}
\nu \int_{(0,T) \times \Gamma_{c\nu} } | \nabla u |^{2} dx dt \rightarrow  0 ,
\\ \label{KatoConditionWeakk}
 u(t,\cdot)  \rightharpoonup u^E (t,\cdot) \text{ in }   \mathcal{H}-w  , \text{ for any }  t \in [0,T ]  .
\end{eqnarray}
\end{Theorem}
Before to start the proof of Theorem \ref{KatoBody}, let us give a few comments and open questions. 

First as mentioned previously, a similar result can be obtained in two dimensions. The proof is even actually simpler. 
Still let us mention that in two dimensions the assumption that the energy is finite is rather restrictive, at least for what concerns the Euler equation, see \cite{GS} for a wider setting. 
Therefore it is natural to wonder whether or not the analysis performed here can be extended to this more general setting. 
In particular it could be that, even under Kato's condition, one misses some interesting dynamics of the Euler case, as for instance the one obtained in the particle limit in \cite{shrinking}, by using the Navier-Stokes equations.

Another natural issue is to extend Theorem \ref{KatoBody} to the case where there are several bodies, or to the case where the fluid-body system occupies a fixed bounded domain. 
This raises some extra technical difficulties as the change of variable performed in Section $2$ does not lead to a time-independent domain. 
Let us also stress that the collision issues can be very different depending on whether one considers the Euler equations or the  Navier-Stokes equations. Let us refer here to \cite{DGVH,hm} and to the references therein.

Also another interesting question raised by Theorem \ref{KatoBody} is about the convergence of the time derivatives of the body's velocity. In particular it was shown in \cite{ogfstt,GS} that  in the Euler case,  the body's velocity is actually analytic in time,  if its boundary is analytic. It is therefore natural to  wonder whether or not the  time derivatives of the  body's velocity for smooth solutions of the Navier-Stokes case also converge to the ones of the Euler case under a Kato type condition.

 It is also probably possible to extend some of the variants of Kato's argument mentioned in the introduction in this setting of a moving body.

\section{Beginning of the proof of Theorem \ref{KatoBody}}

\subsection{Easy part}

As in Kato's original statement, the proof of the necessity of the condition  \eqref{KatoCondition} to get \eqref{ConV}
 is quite easy: if \eqref{ConV} holds true when  $\nu \rightarrow  0$ then it suffices to combine \eqref{NSBodyWeakEnergy}, \eqref{EulerBodyStrongEnergy} and \eqref{cvDonnees} to get that 
\begin{eqnarray}
 \label{KatoConditionTot}
\nu \int_{(0,T) \times \R^2} | D(u) |^{2} dx dt \rightarrow  0 ,
\end{eqnarray}
when  $\nu \rightarrow  0$.
Of course  \eqref{KatoConditionTot} implies \eqref{KatoCondition}.

We obtain similarly that \eqref{ConV}  implies \eqref{KatoConditionVort} and \eqref{KatoConditionGrad} using
Remark \ref{NSBodyWeakRemarkVort} and  Remark \ref{NSBodyWeakRemarkGrad}.

Since it is straightforward that  \eqref{ConV}  implies \eqref{KatoConditionWeakk}, it remains now to see the converse statements. 

Actually let us see that \eqref{KatoConditionWeakk} implies \eqref{KatoCondition} so that it will only remain to prove that either \eqref{KatoCondition}, or \eqref{KatoConditionVort} or  \eqref{KatoConditionGrad} implies \eqref{ConV}.

Thanks to  \eqref{NSBodyWeakEnergy}, we have 
for any $t \in [0,T ]$, using  \eqref{cvDonnees},
\begin{eqnarray*}
  2 \limsup \nu \int_{(0,t) \times \R^{3} } | D(u) |^{2} dx dt 
&\leq& \frac{1}{2}   \| u^{E}_{0}  \|^{2}_{\mathcal{H}} - \liminf \frac{1}{2} \| u (t, \cdot)  \|^{2}_{\mathcal{H}}+ \limsup  \int_{0}^{t}   f_s [u,u]  ds    
\\ &\leq& \frac{1}{2}   \| u^{E}_{0}  \|^{2}_{\mathcal{H}} - \frac{1}{2} \| u^E (t, \cdot)  \|^{2}_{\mathcal{H}}+ \int_{0}^{t}   f_s [u^E ,u^E]  ds  , 
\end{eqnarray*}
using \eqref{KatoConditionWeakk} and Fatou's lemma.
It remains to use \eqref{EulerBodyStrongEnergy} to see that the right hand side above is $0$, what yields \eqref{KatoCondition}.

We will detail how to prove that  \eqref{KatoCondition} implies \eqref{ConV} and then we will explain what modifications lead to the other cases.
We first adapt the construction of a Kato type  ``fake'' layer.

\subsection{A Kato type ``fake'' layer}
\label{Fake}
The goal of this section is to prove the following result, where we make use of the Landau notations $o(1)$ and $O(1)$ for quantities respectively converging to $0$ and bounded with respect to the limit $\nu \rightarrow 0^{+}$.
\begin{Proposition}
\label{PropositionFake}
Under the assumptions of Theorem  \ref{KatoBody}
there exists  $v_{F}\in  C ( [0,T ]  ;  \mathcal{H})$, supported in  $ \Gamma_{c\nu} $, such that 
\begin{eqnarray}
\label{Fake0}
v_{F} = O( 1 ) \text{ in }C ( [0,T ]  \times \R^{3} ),
\\ \label{Fake1}
v_{F} = O( \nu^{\frac{1}{2}}) \text{ in }C ( [0,T ]  ;  \mathcal{H}),
\\ \label{Fake2}
\partial_{t} v_{F} = O( \nu^{\frac{1}{2}}) \text{ in }C ( [0,T ]  ;  \mathcal{H})
\\   \label{Fake4}
\| \nabla v_{F} \|_{ L^{\infty} ([0,T ] ;   L^{2} ( \Gamma_{c\nu} )) } = O( \nu^{-\frac{1}{2}}  ) ,
\\   \label{FakeImp}
d(x) v_{F} = O( \nu ) \text{ in }   L^{\infty} ( [0,T ]  \times \R^{3} ),
\\  \label{Fake7}
u^E  - v_{F}  \in C ( [0,T ]  ;  \mathcal{H}) \cap  L^2 ( [0,T ] ;  \mathcal{V} )
\end{eqnarray}
\end{Proposition}
\begin{proof}
According to \cite{Tosio}, Lemma A1, we get that there exists an antisymmetric 2-tensor field $a_F (t,x)$  on $ [0,T ] \times \R^{3}$  such that, 
\begin{eqnarray}
\label{ajou1}
\div a_F = u^E - u^E_{\mathcal{S}}  \text{ and } a_{F} = 0   \text{ on }  \partial   \mathcal{S}_0 .
\end{eqnarray}
Let us recall that for a smooth  antisymmetric 2-tensor $a$, $\div  a$ denotes the vector field  $\div  a := (\sum_k \partial_k a_{jk} )_{k}$.

Now we introduce a smooth cut-off function $\xi: [ 0,+\infty) \rightarrow [ 0,+\infty) $ such that $\xi(0) = 1$ and  $\xi(r) = 0$ for $r \geq 1$. 
We define $z(x) := \xi ( \frac{ d(x)}{c \nu} )$ and $v_{F}$ by 
\begin{eqnarray}
\label{ajou2}
v_{F} := \div ( z a_F   )   \text{ in }  \mathcal{F}_0 \text{ and } v_{F} := 0  \text{ in } \mathcal{S}_0     . 
\end{eqnarray}

In order to verify that $v_{F}$ satisfies the desired properties, let us  introduce 
$ a_F^{\flat} (t,x) := \frac{1}{d(x)} a_{F} (t,x)$,  $\tilde{\xi}(r) := r  \xi' (r)$ and $\tilde{z}(x) := \tilde{\xi} ( \frac{ d(x)}{c \nu} )$. 
Then, in $\mathcal{F}_0$, 
\begin{eqnarray}
 \label{FakeMother}
v_{F} = z \div a_F +   \tilde{z} a_F^{\flat}  \nabla d.
\end{eqnarray}
First  since $z$ and $\tilde{z}$ are  supported in  $ \Gamma_{c\nu} $ so is $v_{F} $.
Furthermore, using \eqref{ajou1} and that, for $x \in \partial \mathcal{S}_0$,  $z (x)=1 $ and $\tilde{z} (x)= 0$, we get 
\begin{eqnarray}
\label{ajou289}
v_{F} |_{\mathcal{F}_0} = u^E - u^E_{\mathcal{S}}     \text{ on }  \partial  \mathcal{S}_0  . 
\end{eqnarray}

We observe that for any smooth  antisymmetric 2-tensor $a$ the vector field  $\div  a$ is divergence free, as
 $\div \div  a = \sum_j \sum_k \partial_j a_{jk} = 0$.
 Therefore we obtain  that  $v_{F}\in  C ( [0,T ]  ;  \mathcal{H})$.

Moreover $u^E  - v_{F} $ is $H^1$ in $ \mathcal{F}_0$ and in  $\mathcal{S}_0$.
Using again \eqref{ajou289} we get that  $u^E  - v_{F} $ is continuous across $\partial \mathcal{S}_0$.  Therefore it belongs to $L^2 ( [0,T ] ;  \mathcal{V} )$.

The other estimates follow easily from \eqref{FakeMother} if one observes that the functions $z$ and $\tilde{z}$ satisfy the required estimates and that, according to  \eqref{FakeMother}, $v_{F}$ is a slow modulation  (with respect to $\nu$) of  $z$ and $\tilde{z}$ by some  regular functions.
\end{proof}

\section{Core of the proof of Theorem \ref{KatoBody}}
In this section we prove that  \eqref{KatoCondition} implies \eqref{ConV}.
Let us give a few words of caution before entering in the proof: 
\begin{enumerate}
\item We will use the same notation $C$ for various constants (which may change from line to line).
\item For some functions $\phi$ and $\psi$  depending on $(t,x)$, such that for any $t$, $\phi (t,\cdot)$ and $\psi(t,\cdot)$ are in $\mathcal{H}$, we will denote $ (\phi ,\psi  )_{\mathcal{H}}  (t) $ for $(\phi (t,\cdot),\psi  (t,\cdot))_{\mathcal{H}}$.
\item The identities  \eqref{P1} and \eqref{P2} are also true for an unbounded domain, for instance if one substitutes the domain $\mathcal{F}_{0}$ to the domain $\Omega$.
\end{enumerate}
For any $t \in [0,T ]$, we have, thanks to \eqref{NSBodyWeakEnergy}, \eqref{EulerBodyStrongEnergy},  the Cauchy-Schwarz inequality,  \eqref{Fake1} and \eqref{cvDonnees},
\begin{eqnarray}
\nonumber \|  u(t,\cdot) -  u^E  (t, \cdot)   \|_{\mathcal{H}}^{2} &=& \|  u(t,\cdot) \|_{\mathcal{H}}^{2}   +  \|    u^E  (t, \cdot)   \|_{\mathcal{H}}^{2} - 2 (u  ,   u^E   )_{\mathcal{H}} (t)
\\ \nonumber  &\leq& \| u_{0}  \|_{\mathcal{H}}^{2} + \|  u^E_{0}  \|_{\mathcal{H}}^{2}  + 2 \int_{0}^{t}  ( f_s [u^E ,u^E ]  +  f_s [u,u]  ) ds - 2 (u ,   u^E   )_{\mathcal{H}}  (t)
\\  \label{square} &\leq&  2 \|  u^E_{0}  \|_{\mathcal{H}}^{2} + 2 \int_{0}^{t}  ( f_s [u^E ,u^E ]  +  f_s [u,u]  )   ds  - 2 (u ,   u^E   - v_{F} )_{\mathcal{H}}  (t)  + o(1) .
\end{eqnarray}

We now apply \eqref{WeakNS}  to $v = u^E  - v_{F}  $ (what is licit according to  \eqref{Fake7}) to get
\begin{eqnarray*}
(u,  u^E  - v_{F} )_{\mathcal{H}}  (t) -  (u_{0},  u^E_{0}  - v_{F}  |_{t=0})_{\mathcal{H}}   = 
 \int_{0}^{t} \Big[ (u, \partial_{t}  ( u^E  - v_{F} ))_{\mathcal{H}}     + 
  b(u,u,u^E  - v_{F})  
\\ - 2 \nu     \int_{\mathcal{F}_0 }   D(u) :  D(  u^E  - v_{F} )   dx + f_s [ u, u^E ]  \Big] ds .
\end{eqnarray*}
Let us stress that we used above that $f_s [ u, v_{F}] = 0$.
Now using 
\eqref{cvDonnees},  \eqref{NSBodyWeakEnergy}, \eqref{Fake1}, the Cauchy-Schwarz inequality and \eqref{Fake2}  we deduce that 
\begin{eqnarray}
\label{WeakNSappl2}
-2(u,  u^E  - v_{F} )_{\mathcal{H}}  (t)  + 2  \|  u^E_{0} \|^{2}_{\mathcal{H}} = o(1) 
  -2 \int_{0}^{t} \Big[   R(s)   + (u, \partial_{t}   u^E )_{\mathcal{H}}     +  b(u,u,u^E ) + f_s [ u,u^E ]     \Big] ds ,
\end{eqnarray}
where $R$ denotes the time-dependent function:
\begin{eqnarray*}
R &:=&    b(u,u, v_{F} ) + 2 \nu  \int_{\mathcal{F}_0 }D(u) :  D(  u^E  - v_{F} ) dx 
\\ &=& \int_{\mathcal{F}_0 } \Big(     [   (u-  u_\mathcal{S} )  \cdot\nabla v_F ] \cdot u   -  \det (r_u , u , v_F ) \Big) + 2 \nu  \int_{\mathcal{F}_0 }D(u) :  D(  u^E  - v_{F} ) dx .
\end{eqnarray*}
On the other hand we have, for any $t \in [0,T ]$,
\begin{eqnarray*}
  (\partial_{t}  u^E ,  u)_{\mathcal{H}}  =  - b(u^E,u,u^E ) +  f_t [ u^E ,u]   .
\end{eqnarray*}
To see that,  multiply  \eqref{chEuler1} by $v = u$ and integrate by parts in space using \eqref{chEuler1}-\eqref{chEulerSolide2}.

Combining with  \eqref{WeakNSappl2}  we obtain
\begin{eqnarray}
\label{voili}
-2 (u,  u^E  - v_{F} )_{\mathcal{H}}  (t)   + 2  \|  u^E_{0} \|^{2}_{\mathcal{H}} 
= o(1)  - 2 \int_{0}^{t}  \Big[   R(s)   +  b(u-u^E,u,u^E ) +  f_s [u  ,u^E]  +  f_s [ u^E ,u]     \Big] ds 
\end{eqnarray}
Using the property \eqref{Tonga} we get 
\begin{eqnarray*}
-2 (u,  u^E  - v_{F} )_{\mathcal{H}}  (t)   + 2  \|  u^E_{0} \|^{2}_{\mathcal{H}} 
= o(1)  - 2 \int_{0}^{t}   \Big[   R(s)   +  b(u-u^E,u-u^E,u^E ) +   f_s [u  ,u^E]  +  f_s [ u^E ,u]    \Big] ds   ,
\end{eqnarray*}
and, then, using that $(\nabla u^E ) |_{ [0,T ]  \times {\mathcal{F}_0}} \in    C ( [0,T ]  ;  L^2 ({\mathcal{F}_0} , (1+ | x |^2 )^\frac12  dx ) \cap L^{\infty} (\mathcal{F}_0 ))$, we get 
\begin{eqnarray*}
-2 (u,  u^E  - v_{F} )_{\mathcal{H}}  (t)   + 2  \|  u^E_{0} \|^{2}_{\mathcal{H}} 
\leqslant  o(1)  -2 \int_{0}^{t}  R(s) ds 
  + C   \int_{0}^{t}     \| u - u^E \|^{2} _{\mathcal{H}} (s) ds -2  \int_{0}^{t} \Big[  f_s [u  ,u^E]  +  f_s [ u^E ,u]  \Big] ds .
\end{eqnarray*}

Now combining this with \eqref{square} yields
\begin{eqnarray*}
\|  u(t,\cdot) -  u^E  (t, \cdot)   \|_{\mathcal{H}}^{2} &\leq& o(1)  
 -2 \int_{0}^{t} R(s) ds 
 +   C    \int_{0}^{t} 
  \| u - u^E \|^{2}_{\mathcal{H}} (s) ds 
  + 2 \int_{0}^{t}  ( f_s [u^E ,u^E -u ]  +  f_s [u,u- u^E] )   ds  .
  \end{eqnarray*}
  Moreover, combining \eqref{exp1} and \eqref{exp2}, and using again the bounds given by \eqref{NSBodyWeakEnergy} and  \eqref{EulerBodyStrongEnergy}, we obtain, for any $s \in [0,t ]$,
\begin{eqnarray*}
|    f_s [u^E ,u^E -u ]  +  f_s [u,u- u^E]  | \leq C   \| u - u^E \|_{\mathcal{H}} (s)    \sup_{ 0\leq \tilde{s} \leq s}  \| u - u^E \|_{\mathcal{H}} (\tilde{s} ) 
\end{eqnarray*}
As a consequence in order to achieve this part of the proof of Theorem \ref{KatoBody} it only suffices to prove that 
\begin{eqnarray}
\label{voila}
 \int_{0}^{t} R(s) ds  \rightarrow 0 \text{ when } \nu \rightarrow 0 .
 \end{eqnarray}
  In order to prove \eqref{voila} we first decompose $R(t)$ into 
\begin{eqnarray*}
R(t) =  R_1 (t)  + ... +  R_5 (t) ,
\end{eqnarray*}
where 
\begin{eqnarray*}
 R_1   &:=&    - \int  (u -  u_{\mathcal{S}} ) \cdot [   (u-u_{\mathcal{S}}   )  \cdot\nabla  v_{F}  ] dx  ,
  \\ R_2   &:=&   - \int  u_{\mathcal{S}}  \cdot [  (u- u_{\mathcal{S}} )  \cdot\nabla  v_{F}  ] dx  ,
  \\ R_3   &:=& 2 \nu  \int   D(u) :  D(  u^E  ) dx ,
    \\ R_4   &:=& - 2 \nu  \int D(u) :  D(  v_{F} ) dx ,
    \\ R_5    &:=& - \int_{\mathcal{F}_0 }     \det (r_u , u , v_F ) .
\end{eqnarray*}
 Let us emphasize that the integrals in the expressions above, except the one corresponding to $R_{3}$, can be taken over $\Gamma_{c\nu}$, since the fake layer $v_{F} $ is supported in $\Gamma_{c\nu}$.
 In particular we do not have to worry too much about the nondecreasing at infinity of the vector field $  u_{\mathcal{S}}$. 
 However let us gain in comfort by introducing a smooth cut-off function $\chi$ defined on $\mathcal{F}_0$ such that $\chi = 1$ in $\Gamma_{c}$ and  $\chi = 0$ in $\mathcal{F}_0 \setminus \Gamma_{2 c}$.
Let us denote
\begin{eqnarray*}
    \psi_{\mathcal{S}} (t,x) :=  -\frac12 (  \ell (t) \wedge  x  + \frac{1}{2} r  (t)  | x |^{2} )   \text{ and } \tilde{u}_{\mathcal{S}} := \curl ( \chi  \psi_{\mathcal{S}} ) , 
\end{eqnarray*}
and observe that 
\begin{eqnarray*}
 R_1 (t)  =   - \int  (u -   \tilde{u}_{\mathcal{S}} ) \cdot [   (u-  \tilde{u}_{\mathcal{S}}   )  \cdot\nabla  v_{F}  ] dx  \text{ and }
  R_2 (t)  =   - \int   \tilde{u}_{\mathcal{S}}  \cdot [  (u-  \tilde{u}_{\mathcal{S}} )  \cdot\nabla  v_{F}  ] dx  ,
 \end{eqnarray*}
 since  $ \tilde{u}_{\mathcal{S}} = u_{\mathcal{S}}$ on the support of $v_{F}$ for $\nu\leq 1$.
  Moreover, $ \tilde{u}_{\mathcal{S}}$ is a $H^1$  divergence free vector field on $\mathcal{F}_0$ and, using \eqref{NSBodyWeakEnergy}, we  have that
  \begin{eqnarray}
  \label{cool}
\|    \tilde{u}_{\mathcal{S}} \|_{ L^{\infty} ([0,T ] ;   H^{1} ({\mathcal{F}}_0 )) }= O( 1) .
\end{eqnarray}
  Regarding $R_1 (t) $ we first integrate by parts to get
\begin{eqnarray*}
R_1 (t)  &=&    \int  v_{F} \cdot [   (u-\tilde{u}_{\mathcal{S}}   )  \cdot\nabla   (u -  \tilde{u}_{\mathcal{S}} ) ] dx  .
\end{eqnarray*}
Then  we can use the equality \eqref{P2} to obtain 
\begin{eqnarray*}
R_1 (t)  =  2  \int  v_{F}  \cdot     \Big( D(u -  \tilde{u}_{\mathcal{S}} ) (u -  \tilde{u}_{\mathcal{S}} ) \Big) dx  
= 2  \int  v_{F}  \cdot     \Big( D(u ) (u -  \tilde{u}_{\mathcal{S}} ) \Big) dx ,
\end{eqnarray*}
since  $ \tilde{u}_{\mathcal{S}}$ is a rigid velocity on the support of $v_{F}$.

Then 
\begin{eqnarray*}
R_1 (t)  &=& 2  \int  d(x) v_{F}  \cdot     \Big( D(u ) \tau  \Big) dx ,
\end{eqnarray*}
where 
\begin{eqnarray*}
\tau (t,x) := d(x)^{-1 } (u (t,x)  -  \tilde{u}_{\mathcal{S}}  (t,x)  ) .
\end{eqnarray*}
Since  the vector field $u - \tilde{u}_{\mathcal{S}}$ is vanishing on  $\partial \mathcal{S}_0$, according to Hardy's inequality we have,  uniformly in  $t$,
\begin{eqnarray}
 \label{hardy}
\|  \tau    \|_{L^{2} (  \Gamma_{c\nu}  )} &\leq & C \|  \nabla ( u - \tilde{u}_{\mathcal{S}} )  \|_{L^{2} ( \Gamma_{c\nu}   )}  .
\end{eqnarray}
Thus  
\begin{eqnarray*}
 |  R_1 (t)  |   &\leqslant &  C \nu   \|  D(u)  \|_{L^{2} (   \Gamma_{c\nu} ) }  \|  \nabla ( u - \tilde{u}_{\mathcal{S}} )   \|_{L^{2} (\mathcal{F}_0   )}   ,
\end{eqnarray*}
thanks to the Cauchy-Schwarz inequality, \eqref{hardy} and \eqref{FakeImp}. 
Using again the Cauchy-Schwarz inequality with respect to the time integration, that
\begin{eqnarray}
\|  \nabla  ( u - \tilde{u}_{\mathcal{S}} )    \|_{L^{2} (\mathcal{F}_0   )}  = 2 \|  D  ( u - \tilde{u}_{\mathcal{S}} )    \|_{L^{2} (\mathcal{F}_0   )} 
\end{eqnarray}
according to the identity  \eqref{P1}, \eqref{NSBodyWeakEnergy}, \eqref{cool} and \eqref{KatoCondition}, we obtain
\begin{eqnarray}
 \label{R1}
 \int_{0}^{t} |  R_1 (t)  |   ds  \rightarrow 0 \text{ when } \nu \rightarrow 0 .
 \end{eqnarray}
Similarly, we  integrate by parts $ R_2 (t) $ to get 
\begin{eqnarray*}
R_2 (t)  =    \int  v_{F}  \cdot [  (u- \tilde{u}_{\mathcal{S}} )  \cdot\nabla  \tilde{u}_{\mathcal{S}}   ] dx  =   \int  v_{F}  \cdot   [ r(t) \wedge  (u- \tilde{u}_{\mathcal{S}} ) ] dx  ,
\end{eqnarray*}
using that, on the support of $v_{F}$, $\tilde{u}_{\mathcal{S}}$ is given by the formula \eqref{LastReading}.
Then
\begin{eqnarray}
  \label{R2}
 \int_{0}^{t} R_2 (s) ds = O( \nu^{1/2} )  ,
\end{eqnarray}
thanks to \eqref{NSBodyWeakEnergy} and \eqref{Fake1}.

It remains to deal with $R_3$ and  $R_4$.
Using the Cauchy-Schwarz inequality and that  $ u^E  \in L^{\infty} ((0,T) ;  H^1 (\mathcal{F}_0 ) )$, we get
\begin{eqnarray}
\nonumber
| \int_{0}^{t} R_3 (s) ds  |    \leqslant   C \int_{0}^{t}  \nu  \|  D(u)  (s, \cdot)  \|_{L^{2} ( \mathcal{F}_0  )} ds  
   \leqslant   C {t}^{\frac{1}{2}}  \nu  \|  D(u)    \|_{L^{2} ( (0,t) \times \mathcal{F}_0  )} 
\end{eqnarray}
by using again the Cauchy-Schwarz inequality.
Thanks to \eqref{NSBodyWeakEnergy} we obtain
\begin{eqnarray}
  \label{R3}
 | \int_{0}^{t} R_3 (s) ds  |   & \leqslant &  C {t  \nu}^{\frac{1}{2}} .
\end{eqnarray}

Regarding $R_4 (t)$, we have, using \eqref{Fake4}, that
\begin{eqnarray}
  \label{R4}
| \int_{0}^{t} R_4 (s) ds  |  &  \leqslant & C\nu^{\frac{1}{2}}  \|  D(u)  \|_{L^{2} ( (0,t) \times  \Gamma_{c\nu} ) }  = o(1) ,
\end{eqnarray}
thanks to \eqref{KatoCondition}.

Finally, thanks to \eqref{NSBodyWeakEnergy} and \eqref{Fake1} we obtain
\begin{eqnarray}
  \label{R5}
 \int_{0}^{t} R_5 (s) ds    &  =& O(\nu^{\frac{1}{2}} ) .
\end{eqnarray}

Gathering \eqref{R1}-\eqref{R5} we obtain \eqref{voila} and the proof is over.

\section{End of the proof of Theorem \ref{KatoBody}}

In this section we explain how to modify the proof of the previous section in order to obtain that 
 either  \eqref{KatoConditionVort} or  \eqref{KatoConditionGrad} implies \eqref{ConV}.
Of course the idea is to use the weak formulations 
\eqref{WeakNSvort} and \eqref{WeakNSgrad} instead of 
  \eqref{WeakNS} 
and the energy inequalities 
 \eqref{NSBodyWeakEnergyvort} and 
\eqref{NSBodyWeakEnergygrad} instead of  \eqref{NSBodyWeakEnergy}.
Then things go as previously till the treatment of the term $R(t)$ for which we simply  use the identities \eqref{P1curl} and  \eqref{P1grad} instead of \eqref{P1}.

\section{Appendix.}

\subsection{Proof of  Proposition \ref{WeakStrong}}

First observe  that the result of Proposition \ref{WeakStrong} will follow, by an integration by parts in time, from the following claim: for any  $v \in \mathcal{H} \cap C^{\infty }_{c} (\R^{3} )$, for any $t \in [0,T]$, 
\begin{eqnarray}
\label{WeakNSDos}
  (\partial_{t} u, v)_{\mathcal{H}} 
   =  b(u,u,v)    - 2 \nu   \int_{\mathcal{F}_0 }     D(u) :  D(v)  dx + f_t [u,v]  .
\end{eqnarray}

Then we multiply the equation \eqref{chNS1} by $v$ and integrate over $\mathcal{F}_0$:
\begin{eqnarray*}
 \int_{ \mathcal{F}_0}  \frac{\partial u}{\partial t} \cdot v +  \int_{ \mathcal{F}_0}  [ \Big( (u- u_\mathcal{S}    ) \cdot\nabla \Big)u  ]\cdot v
 +  \int_{ \mathcal{F}_0}  ( r(t) \wedge u  )  \cdot v   +  \int_{ \mathcal{F}_0}  \nabla p  \cdot v   
    =  \int_{ \mathcal{F}_0} \nu \Delta u \cdot v  
  + \int_{ \mathcal{F}_0} Q(t)^T \,  g  \cdot v  .
\end{eqnarray*}
We then use some integrations by parts, taking into account \eqref{chNS2} and \eqref{chNS3}, to get
\begin{eqnarray*}
  \int_{ \mathcal{F}_0} [ \Big( (u- u_\mathcal{S}    ) \cdot\nabla \Big)u  ]\cdot v  &=& - \int_{ \mathcal{F}_0} u \cdot  \Big( (u-  u_\mathcal{S}  ) \cdot\nabla \Big)v  ,
 \\  \int_{ \mathcal{F}_0}   ( r(t) \wedge u  )  \cdot v &=&   \int_{ \mathcal{F}_0}  \det (r,u,v)   ,
 \\  \int_{ \mathcal{F}_0}  \nabla p  \cdot v &=&   \int_{ \partial  \mathcal{S}_0}  p  n \cdot v  ,
 \\   \int_{ \mathcal{F}_0} \nu \Delta u \cdot v &=&  2\nu  \int_{ \partial \mathcal{S}_0}  (D(u)  v) \cdot  n - 2 \nu \int_{ \mathcal{F}_0} D(u) :  D(v) ,
 \\  &=&  2\nu  \int_{ \partial \mathcal{S}_0}  (D(u)  n) \cdot v - 2 \nu \int_{ \mathcal{F}_0} D(u) :  D(v) ,
\end{eqnarray*}
since $D(u) $ is symmetric.
Then we observe that 
\begin{eqnarray*}
 \int_{ \partial  \mathcal{S}_0}  p  n \cdot v - 2\nu  \int_{ \partial  \mathcal{S}_0}   (D(u)  n) \cdot v
 &=& - \ell_{v} \cdot \int_{ \partial \mathcal{S}_0 } \sigma n \, ds   - r_{v}  \cdot  \int_{ \partial \mathcal{S}_0 }  x \wedge \sigma n \, ds
\\ \label{rhs}
&=& m  \ell_{v}  \cdot   \ell'   + \mathcal{J}_0    r_{v}  \cdot r'  -  \det (ml,r, \ell_{v})  -  \det (  \mathcal{J}_0  r ,r,    r_{v}) - m  \ell_{v}  \cdot  Q(t)^T \,  g  ,
\end{eqnarray*}
thanks to \eqref{chSolide1}-\eqref{chSolide2}. 

Finally we have the following simplification of the gravity contribution, what corresponds to the Archimedes' principle, 
\begin{eqnarray*}
 \int_{ \mathcal{F}_0}  Q(t)^T \,  g   \cdot v &=&  \int_{ \R^2}  Q(t)^T \,  g   \cdot v -  \int_{ \mathcal{S}_0}  Q(t)^T \,  g   \cdot v  
 \\ &=&  \int_{ \R^2}  \nabla (Q(t)^T \,  g \cdot x)   \cdot v -  \int_{ \mathcal{S}_0}  Q(t)^T \,  g   \cdot v  
  \\ &=& -  \int_{ \mathcal{S}_0}  Q(t)^T \,  g   \cdot v  ,
\end{eqnarray*}
  since $v$ is divergence free. Moreover in $ \mathcal{S}_0$, $v =  \ell_{v} + r_{v} \wedge x $ so that 
\begin{eqnarray*}
 \int_{ \mathcal{F}_0}  Q(t)^T \,  g   \cdot v
  &=& -   Vol(\mathcal{S}_0) Q(t)^T \,  g   \cdot  ( \ell_{v} + r_{v} \wedge x_0 )  ,
\end{eqnarray*}
by definition of $x_0$.

Gathering all these equalities yields \eqref{WeakNSDos}.

\subsection{Proof of Remark \ref{NSBodyWeakRemarkVort} and of Remark \ref{NSBodyWeakRemarkGrad}}

Let us now explain briefly how to modify the previous calculations in order to prove the claims in Remark \ref{NSBodyWeakRemarkVort}  and in Remark \ref{NSBodyWeakRemarkGrad}.
First, for any  $v \in \mathcal{H} \cap C^{\infty }_{c} (\R^{3} )$, for any $t \in [0,T]$, 
\begin{eqnarray*}
 \int_{ \mathcal{F}_0} \nu \Delta u \cdot v &=&  - \nu  \int_{ \mathcal{F}_0}  v \cdot (\nabla \wedge \omega ) ,
 \end{eqnarray*}
 where $\omega := \nabla \wedge u$ is the vector in $\R^3$ canonically associated to the $3 \times 3$ matrix  $\curl  u$. 
 Now, using the following formula, for two smooth enough vector fields $a$ and $b$
 \begin{eqnarray}
 \label{divrot}
 - a  \cdot (\nabla \wedge b) =  \div (a \wedge b ) - b  \cdot (\nabla \wedge a) ,
\end{eqnarray}
we get
 \begin{eqnarray*}
 \int_{ \mathcal{F}_0} \nu \Delta u \cdot v &=&  \nu  \int_{ \partial \mathcal{S}_0} (v \wedge \omega )  \cdot n  - \nu  \int_{ \mathcal{F}_0}  \omega   \cdot (\nabla \wedge v) 
 \\ &=&  \nu  \int_{ \partial \mathcal{S}_0} (v \wedge \omega )  \cdot n  - 2 \nu  \int_{ \mathcal{F}_0}  \curl u  : \curl v  .
 \end{eqnarray*}
On the other hand, we have classically
\begin{eqnarray*}
 \int_{ \mathcal{F}_0} \nu \Delta u \cdot v &=&-  \int_{ \mathcal{F}_0}  \nabla u : \nabla v  +   \nu  \int_{ \partial \mathcal{S}_0} v (\nabla u )^T  n .
\end{eqnarray*}
Therefore, in order to prove that the weak formulation can be modified as stated in  Remark \ref{NSBodyWeakRemarkVort} and in Remark \ref{NSBodyWeakRemarkGrad}, it is sufficient to prove 
\begin{eqnarray}
\label{toprove}
  \int_{ \partial \mathcal{S}_0} (v \wedge \omega )  \cdot n =  \int_{ \partial \mathcal{S}_0}  v (\nabla u )^T  n  = 2 \int_{ \partial \mathcal{S}_0}   (D(u)  n ) \cdot v .
\end{eqnarray}
Thus, let us write:
 \begin{eqnarray}
 \label{yep1}
  \int_{ \partial \mathcal{S}_0} (v \wedge \omega )  \cdot n =  \int_{ \partial \mathcal{S}_0}  v \cdot ( \omega \wedge n ) = \ell_{v}  \cdot   \int_{ \partial \mathcal{S}_0}   \omega \wedge n  +  r_{v}\cdot   \int_{ \partial \mathcal{S}_0}  x \wedge (\omega \wedge n  )  ,
 \\   \label{yep2}
 \int_{ \partial \mathcal{S}_0} v (\nabla u )^T  n   =  \ell_{v}  \cdot  \int_{ \partial \mathcal{S}_0} (\nabla u )^T  n +  r_{v}\cdot   \int_{ \partial \mathcal{S}_0}  x \wedge ((\nabla u )^T  n) .
  \end{eqnarray}
Now, observe that
\begin{eqnarray*}
\int_{ \partial \mathcal{S}_0} n_j \partial_i v_j = \int_{ \partial \mathcal{S}_0} n \cdot \Big(\nabla \wedge (v \wedge e_i )\Big) =  \int_{  \mathcal{F}_0} \div [ \nabla \wedge (v \wedge e_i )] = 0,
\\ \int_{ \partial \mathcal{S}_0} x \wedge ( n_j \partial_i v_j )= \int_{ \partial \mathcal{S}_0} x \wedge [  n \cdot \Big(\nabla \wedge (v \wedge e_i )\Big) ]
= \int_{ \partial \mathcal{S}_0}  n \cdot [  x \wedge  \Big(\nabla \wedge (v \wedge e_i )\Big) ]
 = \int_{  \mathcal{F}_0} \div [  x \wedge  \Big(\nabla \wedge (v \wedge e_i )\Big) ]
= 0, 
\end{eqnarray*}
by using again \eqref{divrot}.
Thus we get
 \begin{eqnarray*}
\int_{ \partial \mathcal{S}_0}   \omega \wedge n =   \int_{ \partial \mathcal{S}_0}  (\nabla u )^T  n        = 2 \int_{ \partial  \mathcal{S}_0}   D(u)  n ,
\\ \int_{ \partial \mathcal{S}_0}  x \wedge (\omega \wedge n  )  =   \int_{ \partial \mathcal{S}_0}     x \wedge  ( (\nabla u )^T  n  )       = 2 \int_{ \partial \mathcal{S}_0}  x \wedge (  D(u)   n  )   .
\end{eqnarray*}
Combining this with \eqref{yep1} and \eqref{yep2} yields \eqref{toprove}.
\\
\\ {\bf Acknowledgements.} The author was partially supported by the Agence Nationale de la Recherche, Project CISIFS, 
grant ANR-09-BLAN-0213-02. 

\end{document}